\documentclass[11pt]{article}
\usepackage{amsfonts,amsmath,amssymb}
\usepackage{amsthm}
\usepackage{enumerate}
\usepackage[dvipdfmx]{graphicx}
\usepackage{epic}
\usepackage{eepic}
\usepackage{tikz}
\usepackage{multirow}
\usepackage{fouridx}
\usepackage{thmtools,thm-restate}

\usepackage[hypertex]{hyperref}%
\usepackage{comment}

\usepackage{bm}

\makeatletter
    
    \@addtoreset{equation}{section}
  \makeatother

\newcommand{\R }{\mathbb{R}}
\newcommand{\C }{\mathbb{C}}
\newcommand{\N }{\mathbb{N}}
\newcommand{\K }{\mathbb{K}}
\newcommand{\Z }{\mathbb{Z}}
\newcommand{\pr}[1]{\mathbb{P}^{#1}}

\newcommand{\GL}[1]{GL(#1,\K)}

\newcommand{\ba }{\bm{a}}

\newcommand{\bd }{\bm{d}}
\newcommand{\be }{\bm{e}}
\theoremstyle{plain} 
\newtheorem{thm}{Theorem}[section] 
\newtheorem{thmalph}{Theorem}

\newtheorem{fact}[thm]{Fact}

  \newtheorem{prop}[thm]{Proposition} %
  \newtheorem{cor}[thm]{Corollary}%
 
\theoremstyle{definition} 
  \newtheorem{df}[thm]{Definition}
\newtheorem*{notation}{Notation}
  \newtheorem{exmp}[thm]{Example}

  \newtheorem{remark}[thm]{Remark}%

\usepackage[truedimen,margin=30truemm]{geometry}
\allowdisplaybreaks

\title{Description of infinite orbits on multiple projective spaces\footnote{This work was supported by JSPS KAKENHI Grant Number JP16J06813.}}
\author{Naoya SHIMAMOTO\footnote{Email adress: naoyas@ms.u-tokyo.ac.jp}}
\date{}

\begin{document}

\maketitle
\begin{center}
\textit{Graduate School of Mathematical Sciences, University of Tokyo, 3-8-1 Komaba, Meguro-ku, Tokyo 153-8914, Japan}
\end{center}

\begin{abstract}
Let $G$ be the general linear group of the degree $n\geq 2$ over the field $\mathbb{K}=\mathbb{R}$ or $\mathbb{C}$. 
In this article, we give a description of orbit decomposition of the multiple projective space $G^m/P^m$ under the diagonal action of $G$ where $P$ is the maximal parabolic subgroup of $G$ such that $G/P\cong\mathbb{P}^{n-1}\mathbb{K}$. 
We also construct representatives of orbits. 
If $m\geq 4$, the number of orbits is infinite, and we give a description of those uncountably many orbits. 

\bigskip



\end{abstract}

\tableofcontents

\section{Introduction}

The starting point of this research is the orbit decomposition of the flag variety $G/P$ under the action of $H$ where $G$ is a reductive group, $H$ is a closed subgroup of $G$, and $P$ is a parabolic subgroup of $G$. 
It is related to the representation theory. 
For example, for a real reductive algebraic group $G$, its minimal parabolic subgroup $P_G$ (resp.\ a Borel subgroup $B_G$ of the complexification $G_c$ of $G$), and its algebraically defined closed subgroup $H$, Kobayashi-Oshima \cite{ko} proved that the finiteness (resp. boundedness) of the multiplicities of irreducible admissible representations $\pi$ of $G$ in the induced representations $\mathrm{Ind}_H^G(\tau)$ for finite-dimensional irreducible representations $\tau$ of $H$ is equivalent to the existence of open $H$-orbits (resp. $H_c$-orbits) on $G/P_G$ (resp. $G_c/B_G$). 
A homogeneous space $G/H$ satisfying these equivalent condition is called a real spherical variety (resp.\ spherical variety). 
A number of people, for instance, Brion, Vinberg, Kimelfeld, Bien, and Matsuki proved that $G/H$ is real spherical (resp.\ spherical) if and only if the number of $H$-orbits (resp.  $H_c$-orbits) on the flag variety $G/P_G$ (resp.\ $G_c/B_G$) is finite \cite{br1, br2, morb}. 
In general, if there are only finitely many orbits on a variety under the action of an algebraic group, then there exists an open orbit. 
On the other hand, the converse does not hold in general. 
For instance, for a non-minimal parabolic subgroup $P$ of $G$, it may occur that $H$ has infinitely many orbits and some open orbits on $G/P$ simultaneously. 
In the setting of this article: $P$ is the maximal parabolic subgroup of $\GL{n}$ such that $\GL{n}/P\cong\pr{n-1}\K$ where $n\geq 2$ and $\K=\R$ or $\C$, $P$ is not a minimal parabolic subgroup if and only if $n\geq 3$. 
From Corollaries \ref{thm:open} and \ref{thm:finite}, we can see that our setting contains the following 3 situations on $\mathrm{diag}(\GL{n})$-orbit decomposition on the multiple flag variety $\GL{n}^m/P^m$: 
\begin{enumerate}
\item the number of $\mathrm{diag}(\GL{n})$-orbits on $\GL{n}^m/P^m$ is finite (in patricular, there exists an open orbit): $m\leq 3$, 
\item the number of orbits is infinite, but there exists an open orbit: $4\leq m\leq n+1$, 
\item there are no open orbits (in patricular, the number of orbits is infinite): $n+2\leq m$.  
\end{enumerate}
Not only for the minimal parabolic subgroup $P_G$ of $G$, but for a general parabolic subgroup $P$ of $G$, the existence of open $H$-orbits and the finiteness of $H$-orbits on the flag variety $G/P$ have some constraints on the branching law between representations of $H$ and representations of $G$ induced from characters of $P$ \cite{tt}. 

On the other hand, if $(G,H)$ is a symmetric pair, the number of $H$-orbits on a flag variety $G/P$ is always finite. 
Description of the orbit decomposition in this setting was given by Matsuki \cite{mp,mpp}. 
The pair $(G^m, \mathrm{diag}(G))$ we are dealing with in this article is a symmetric pair if and only if $m=2$ by the involutive action of $\mathcal{S}_2$, and the $\mathrm{diag}(G)$-orbit decomposition of $(G\times G)/(P_1\times P_2)$ is described in terms of Bruhat decomposition. 
In fact, $(G^m, \mathrm{diag}(G))$ is no more a symmetric pair if $m\geq 3$, but $\mathrm{diag}(G)$ is regarded as the fixed point subgroup under the action of $\mathcal{S}_m$ naturally defined on the real reductive group $G^m$. 

Kobayashi-Oshima \cite{ko} also proved that the finiteness of the dimension of the space of symmetry breaking operators ($H$-intertwining operators from irreducible admissible representations $\pi$ of $G$ to irreducible admissible representations $\tau$ of $H$) is equivalent to the existence of open $P_H$-orbits on $G/P_G$ where $H$ is reductive and $P_H$ is its minimal parabolic subgroup. 
Furthermore, such pairs $(G,H)$ are classified by Kobayashi-Matsuki \cite{km} under the assumption that $(G,H)$ are symmetric pairs. 

Not only the finiteness of orbits or the existence of open orbits, describing explicit representatives of orbits is also related to the representation theory. 
For the symmetric pair $(G,H)=(O(n+1,1),O(n,1))$ which occurs in the classification given by \cite{km}, the description of $P_H$-orbits on $G/P_G$ plays a role in the explicit construction and the classification of symmetry breaking operators \cite{ks,ks2}. 
Our research comes from these motivations. 

The main target of this article is included in the description of $\mathrm{diag}(G)$-orbits on the multiple flag variety $G^m/(P_1\times P_2\times \cdots\times P_m)$. 
If it satisfies the condition that the number of $\mathrm{diag}(G)$-orbits is finite, we call it to be of finite type. 
Magyar-Weyman-Zelevinsky classified multiple flag varieties of finite type where $G$ is a general linear group \cite{mwzA} and $G$ is a symplectic group \cite{mwzC}. 
Matsuki classified ones of finite type where $G$ is an orthogonal group of the odd degree \cite{mB}, and gave explicit representatives of orbits in some cases \cite{mBex}. 
These results were given by using the notion of splittings of multiple flags into indecomposable ones. 
Multiple flags (elements of multiple flag varieties) admits the notion of direct sums, dimension vectors, and isomorphisms which are compatible with each others. 
In this notion, the double coset $\mathrm{diag}(G)\backslash G^m/(P_1\times P_2\times \cdots\times P_m)$ is identified with the set of all isomorphism classes of multiple flags having a certain dimension vector $\ba$ which is determined by the shape of the parabolic subgroup $P_1\times\cdots\times P_m\subset G^m$. 
Each isomorphism class of multiple flags has a unique splitting into indecomposable ones up to isomorphisms, hence it determines a unique splitting of the dimension vector $\ba$ into smaller ones of indecomposable multiple flags (see Section 3.1). 

Actually, under the condition that the multiple flag variety is of finite type, indecomposable multiple flags having each dimension vector smaller than $\ba$ are always unique up to isomorphisms. 
Hence, describing the double coset $\mathrm{diag}(G)\backslash G^m/(P_1\times P_2\times \cdots\times P_m)$ is naturally identified with the set of all splittings of the dimension vector $\ba$ into smaller ones (see Section 3.2).  
The essential difference of our work from their studies is the lack of this property. 

If we try to describe orbits on multiple flag varieties of infinite type, we cannot use these techniques straightforward. 
In this article, we focus on the orbit decomposition of an $m$-tuple flag variety $G^m/P^m$ under the diagonal action of $G$ where $G=\GL{n}$, $n\geq 2$, and $P$ is the maximal parabolic subgroup such that $G/P\cong\pr{n-1}\K$. 
Here, the double coset $\mathrm{diag}(G)\backslash G^m/P^m$ is identified with the set of all isomorphism classes of multiple flags having the dimension vector $\ba_{n,m}=({}^t(1,n-1)^m)\in M_{2,m}(\N)$. 
In this setting, there are only finitely many orbits if and only if $m\leq 3$, and there exists an open orbit if and only if $n+1\geq m$. 
Hence, this setting can be regarded as the simplest case among the ones where some open orbits and uncountably many orbits exist simultaneously on flag varieties. 
Then we have another problem: for each vector $\bd$ smaller than $\ba_{n,m}$, we have to describe isomorphism classes of indecomposable multiple flags whose dimension vectors are $\bd$, which can exist uncountably many. 
In this article, we describe them by taking a certain open dense $\mathrm{diag}(G)$-stable subset of the multiple flag variety (see Proposition \ref{thm:indrep}).

\begin{thmalph}[see Theorems \ref{thm:main} and \ref{thm:rep}] 
There exists a natural surjection \[\pi\colon\mathrm{diag}(G)\backslash G^m/P^m \twoheadrightarrow \mathcal{P}_{n,m}\]
defined in (\ref{eq:mainsurj}) onto the finite set $\mathcal{P}_{n,m}$ (see Definition \ref{thm:pnm}) with the following property: for each element 
\begin{equation} \label{eq:pnmelement1}
p=\left(\left\{I_a\right\}_{a=1}^{A}, \left\{J_b\right\}_{b=1}^{B}, \left\{ (K_c,r_c) \right\}_{c=1}^{C} \right)
\end{equation}
of $\mathcal{P}_{n,m}$, there exists an open dense embedding \[\prod _{c=1}^C(\pr{r_c-1}\K)^{\# K_c-r_c-1} \hookrightarrow \pi^{-1}(p)\] into the fibre of $p$ with explicit formulae to give representatives of orbits given in (\ref{eq:emb}). 
\end{thmalph} 
Each element of the finite set $\mathcal{P}_{n,m}$ is a tuple as in (\ref{eq:pnmelement1}) where $\{I_a, J_b, K_c\}_{a,b,c=1}^{A,B,C}$ is a partition of the set $\{1,2,\ldots,m\}$ and $\{r_c\}_{c=1}^C$ is a tuple of positive integers satisfying some conditions depending on $n$ (see Definition \ref{thm:pnm}). 
Here, $\mathcal{P}_{n,m}$ is naturally identified with the set of all splittings of the dimension vector $\ba_{n,m}$ (see Proposition \ref{thm:pnmid}). 
Each open dense embedding into the fibre of $\pi$ parametrises indecomposable multiple flags having the same dimension vectors. 
Although project spaces are compact, they are embedded open densely into fibres of $\pi$, since the fibres are not Hausdorff. 
Here is an example for small dimension case: 
\begin{exmp} \label{thm:hausdorff}
If $(n,m)=(2,4)$, then $p_g:=\left(\emptyset,\emptyset,\{(\{1,2,3,4\},2)\}\right)$ is an element of $\mathcal{P}_{2,4}$. 
Then we have an open dense embedding from $\pr{1}\K$ to $\pi^{-1}(p_g)$. 
It is homeomorphic to the quotient space of $\pr{1}\K\amalg\pr{1}\K:=\{1,2\}\times\pr{1}\K$ according to the equivalence relation defined by 
\[(i,x)\sim (j,y)\Leftrightarrow \begin{cases} x=y & \textrm{if }x\neq 0,1,\infty \\ (i,x)=(j,y) & \textrm{if }x=0,1,\infty \end{cases}\]
for $i,j=1,2$ and $x,y\in\pr{1}\K$. 
Hence $\pi^{-1}(p_g)$ is not Hausdorff, and $\pr{1}\K$ is embedded open densely.  
\end{exmp}
\begin{notation}
$\N=\{0,1,2,\ldots\}$. 
\end{notation}

\section{Main results}

Let $G$ be the general linear group of the degree $n\geq 2$ over the field $\K=\R$ or $\C$, and $P$ be its maximal parabolic subgroup such that $G/P\cong\pr{n-1}\K$. 
We consider the multiple flag variety $G^m/P^m$. 
Our main problem is to describe the $\mathrm{diag}(G)$-orbit decomposition of this multiple flag variety. 
We have $G^m/P^m\cong (\pr{n-1}\K)^m$. 
Hence our main problem is equivalent to describing the orbit decomposition of $(\pr{n-1}\K)^m$ under the diagonal action of $\GL{n}$. 
In Section \ref{surjection}, we state the first half of the main result: we give the natural surjection $\pi$ from $\mathrm{diag}(G)\backslash G^m/P^m$ to a finite set $\mathcal{P}_{n,m}$ (see Theorem \ref{thm:main}), and in Section \ref{representatives}, we state the last half of the main result: we give representatives of orbits in each fibre of $\pi$ (see Theorem \ref{thm:rep}). 

\subsection{Decomposition of $G^m/P^m$ into a finite number of $\mathrm{diag}(G)$-stable subsets \label{surjection}}

We begin by introducing a finite set $\mathcal{P}_{n,m}$. 
\begin{df} \label{thm:pnm}
Let $n,\,m$ be positive integers. 
We define $\mathcal{P}_{n,m}$ as the set of all tuples 
\begin{equation} \label{eq:pnmelement}
p=\left(\left\{I_a\right\}_{a=1}^{A}, \left\{J_b\right\}_{b=1}^{B}, \left\{ (K_c,r_c) \right\}_{c=1}^{C} \right)
\end{equation}
satisfying the following conditions: 
\begin{gather}
\left(\coprod_{a=1}^{A}I_{a}\right)\sqcup\left(\coprod_{b=1}^{B}J_{b}\right)\sqcup\left(\coprod_{c=1}^{C}K_{c}\right)=\{1,2,\cdots, m\}, \label{eq:pnm1} \\
1\leq\#I_a,\,3\leq\#J_b,\,4\leq r_c+2\leq\#K_c,\,r_c\in\N\textrm{  for }1\leq a\leq A,\,1\leq b\leq B,\,1\leq c\leq C, \label{eq:pnm2} \\
0\leq r(p)\leq n, \label{eq:pnm3}
\end{gather}
where $r(p):=A+\sum_{b=1}^B(\#J_b-1)+\sum_{c=1}^Cr_c$. 
\end{df}
From now on, we fix $n$ and $m$, and write $\mathcal{P}:=\mathcal{P}_{n,m}$. 
Now, we define a tuple of integers
\begin{equation} \label{eq:dkr}
\bd(K,r):=\left(r;\delta_{1,K},\delta_{2,K},\ldots,\delta_{m,K}\right)
\end{equation}
for $K\subset\{1,2,\ldots,m\}$ and $r\in\N$ where $\delta_{i,K}=1$ if $i\in K$ and $0$ if $i\notin K$ for $1\leq i\leq m$. 
More generally, let $\bd=(r;d^{(1)},d^{(2)},\ldots,d^{(m)})=(r;d^{(i)})_{i=1}^m$ be a tuple of non-negative integers satisfying $r\geq d^{(i)}$ for $1\leq i\leq m$. 
For an $r$-dimensional vector space $V$ (throughout this article, we assume all vector spaces are over $\K$), we define $Fl_{\bd}(V)$ as the direct product of Grassmannians $\prod_{i=1}^m\mathrm{Gr}_{d^{(i)}}(V)$ with elements of the form $\left(V;V^{(1)},V^{(2)}, \ldots, V^{(m)}\right)=(V;V^{(i)})_{i=1}^m$ where $V^{(i)}$ is a $d^{(i)}$-dimensional subspace of $V$ for each $1\leq i\leq m$. 
We write $\dim D=\bd$ and call $D$ a ($\bd$-dimensional) $m$-tuple flag if there exists an $r$-dimensional vector space $V$ and $D\in Fl_{\bd}(V)$. 
Let $\bd=(r;d^{(i)})_{i=1}^m$ and $\be=(s;e^{(i)})_{i=1}^m$ satisfy $r\geq d^{(i)},\,s\geq e^{(i)}$ for $1\leq i\leq m$, then the direct sum of two $m$-tuple flags $D=(V;V^{(i)})_{i=1}^m\in Fl_{\bd}(V)$ and $E=(W;W^{(i)})_{i=1}^m\in Fl_{\be}(W)$ is defined by $D\oplus E:=(V\oplus W;V^{(i)}\oplus W^{(i)})_{i=1}^m\in Fl_{\bd+\be}(V\oplus W)$ where $\bd+\be:=(r+s;d^{(i)}+e^{(i)})_{i=1}^m$. 
Furthermore, if there exists a linear isomorphism $f\colon V\xrightarrow{\sim} W$ such that $f(V^{(i)})=W^{(i)}$ for $1\leq i\leq m$, then we say that $D$ and $E$ are isomorphic. 
We say an $m$-tuple flag $D$ is indecomposable if it is not isomorphic to any direct sum of two non-trivial $m$-tuple flags (we say $D$ to be a trivial $m$-tuple flag if $\dim D=(0;0^m)$). 
In our setting, $G^m/P^m$ is naturally identified with $Fl_{\bd(\{1,2,\ldots,m\},n)}(\K^n)=Fl_{(n;1^m)}(\K^n)$. 

Now, let $F$ be an $m$-tuple flag and $p\in\mathcal{P}$ of the form (\ref{eq:pnmelement}), we say $F$ to be $p$-constructible if $F$ splits into indecomposable $m$-tuple flags as follows: 
\begin{equation} \label{eq:splitspecific}
F\underset{\mathrm{isom}}{\sim}\left(\bigoplus_{k=1}^{n-r(p)}F_{\emptyset,k}\right)\oplus\left(\bigoplus_{a=1}^A F_a\right)\oplus\left(\bigoplus_{b=1}^B D_b\right)\oplus\left(\bigoplus_{c=1}^C E_c\right) \textrm{  satisfying the condition }(p).
\end{equation}
The condition ($p$) is defined as follows: for $1\leq k\leq n-r(p)$, $1\leq a\leq A$, $1\leq b\leq B$, $1\leq c\leq C$, 
\begin{equation*}
\begin{cases}
\dim F_{\emptyset,k}=\bd(\emptyset,1),\, \dim F_a=\bd(I_a,1),\, \dim D_b=\bd(J_b,\#J_b-1),\, \textrm{and}\, \dim E_c=\bd(K_c,r_c), & \\
F_{\emptyset,k},\, F_a,\, D_b,\, \mathrm{and}\, E_c\, \textrm{are indecomposable}.  &
\end{cases}
\end{equation*}
With these notations, we claim the first main result: 
\begin{thm}  \label{thm:main}
let $G=\GL{n}$, $P$ be its maximal parabolic subgroup such that $G/P\cong\pr{n-1}\K$, and $\mathcal{P}=\mathcal{P}_{n,m}$ be the finite set introduced in Definition \ref{thm:pnm}. 
For any element $F\in G^m/P^m$, there exists a unique element $p_F\in\mathcal{P}$ such that $F$ is $p_F$-constructible, which induces a surjective map
\begin{align} \label{eq:mainsurj}
\pi\colon \mathrm{diag}(G)\backslash G^m/P^m\rightarrow\mathcal{P},\, \mathrm{diag}(G)\cdot F \mapsto p_F. 
\end{align}
\end{thm}
\begin{remark}
Alternatively, we can also define $\pi$ as the composition of the surjection (\ref{eq:surjabstract}) and the inverse of the bijection (\ref{eq:pnmid}). 
\end{remark}
We give some examples of the set $\mathcal{P}$ and the map $\pi$ for specific $(n,m)$, and compare them with descriptions of the orbit decompositions obtained by direct computations.  
\begin{exmp} \label{thm:23}
If $(n,m)=(2,3)$, we are considering the case where $G=\GL{2}$ acts diagonally on $X_{2,3}:=Fl_{\bd(\{1,2,3\},2)}(\K^2)=Fl_{(2;1,1,1)}(\K^2)\cong(\pr{1}\K)^3$. 
In this case, the elements of $\mathcal{P}=\mathcal{P}_{2,3}$ are
\begin{align}
\begin{split} \label{eq:23ex}
&p_c:=\left(\{\{1,2,3\}\},\emptyset,\emptyset\right), \\
&p_i:=\left(\{\{i\},\{1,2,3\}\setminus\{i\}\},\emptyset,\emptyset\right) \textrm{ for }i=1,2,3, \\
&p_o:=\left(\emptyset,\{\{1,2,3\}\},\emptyset\right). 
\end{split}
\end{align}
Since any three distinct points in $\pr{1}\K$ can be transposed to $\K e_1,\, \K e_2$, and $\K(e_1+e_2)$ simultaneously by the action of $G$ where $\{e_1,e_2\}$ is the standard basis of $\K^2$, we can see that the $\mathrm{diag}(G)$-orbit decomposition of the multiple flag variety $X_{2,3}$ is described as follows: 
\begin{align}
X_{2,3}&={\cal O}_c\amalg{\cal O}_1\amalg{\cal O}_2\amalg{\cal O}_3\amalg{\cal O}_o \textrm{  where} \label{eq:23decomp} \\
{\cal O}_c&:=\left\{(\K^2;\K v_1,\K v_1,\K v_1)\in X_{2,3}\right\}, \notag \\
{\cal O}_1&:=\left\{(\K^2;\K v_1,\K v_2,\K v_2)\in X_{2,3}\left|\K v_1\neq \K v_2\right.\right\} \textrm{  respectively for }i=2,3 \notag \\
{\cal O}_o&:=\left\{(\K^2;\K v_1,\K v_2,\K v_3)\in X_{2,3}\left|\K v_1,\,\K v_2, \textrm{ and }\K v_3\textrm{ are distinct}\right.\right\}. \notag
\end{align}
\begin{enumerate}
\item For $(\K^2;\K v_1,\K v_1,\K v_1)\in{\cal O}_c$, we have an splitting into indecomposable triple flags as
\[(\K^2;\K v_1,\K v_1,\K v_1)=(\K w; 0,0,0)\oplus(\K v_1;\K v_1,\K v_1,\K v_1)\]
by taking $w\in\K^2\setminus\K v_1$, and $\dim (\K w; 0,0,0)=\bd(\emptyset,1)$, $\dim(\K v_1;\K v_1,\K v_1,\K v_1)= \bd(\{1,2,3\},1)$. 
Hence $\pi({\cal O}_c)=p_c$. 
\item For $(\K^2;\K v_1,\K v_2,\K v_2)\in{\cal O}_1$, we have an splitting into indecomposable triple flags as
\[(\K^2;\K v_1,\K v_2,\K v_2)=(\K v_1; \K v_1,0,0)\oplus(\K v_2; 0,\K v_2,\K v_2),\]
and $\dim (\K v_1; \K v_1,0,0)=\bd(\{1\},1)$, $\dim(\K v_2; 0,\K v_2,\K v_2)=\bd(\{2,3\},1)$. 
Hence $\pi({\cal O}_1)=p_1$. 
Similarly, $\pi({\cal O}_i)=p_i$ for $1\leq i\leq 3$. 
\item If $(\K^2;\K v_1,\K v_2,\K v_3)\in{\cal O}_c$, then it is indecomposable, and $\dim (\K^2;\K v_1,\K v_2,\K v_3)=\bd(\{1,2,3\},2)$. 
Hence $\pi({\cal O}_o)=p_o$. 
\end{enumerate}
\end{exmp}
\begin{exmp}[see Example \ref{thm:hausdorff}] \label{thm:24}
If $(n,m)=(2,4)$, we are considering the case where $G=\GL{2}$ acts diagonally on $X_{2,4}:=Fl_{\bd(\{1,2,3,4\},2)}(\K^2)=Fl_{(2;1,1,1,1)}(\K^2)\cong(\pr{1}\K)^4$, and the elements of $\mathcal{P}=\mathcal{P}_{2,4}$ are
\begin{align}
\begin{split} \label{eq:24ex}
p_c&:=\left(\{\{1,2,3,4\}\},\emptyset,\emptyset\right), \\
p_{\{I_1,I_2\}}&:=\left(\{I_1,I_2\},\emptyset,\emptyset\right) \textrm{  for }\emptyset\neq I_1, I_2 \textrm{ such that }I_1\amalg I_2=\{1,2,3,4\}, \\
p_g&:=\left(\emptyset,\emptyset,\left\{\left(\{1,2,3,4\},2\right)\right\}\right). 
\end{split}
\end{align}
Conisider a $\mathrm{diag}(G)$-stable subset
\begin{equation} \label{eq:xgdef}
X_g:=\left\{(\K^2;\K v_1,\K v_2,\K v_3, \K v_4)\in X_{2,4}\left|\#\{\K v_1,\,\K v_2,\, \K v_3,\, \K v_4\}\geq 3\right.\right\}, 
\end{equation}
For three distinct points in $\pr{1}\K$, there exists $p\in\GL{2}$ such that $g\in\GL{2}$ fixes the three points if and only if $p^{-1}gp$ is a scalar matrix. 
Hence we can see that the decomposition of $X_g$ into $\mathrm{diag}(G)$-orbits is described as follows: 
\begin{align}
X_g&=\coprod_{I\subset \{1,2,3,4\},\, \#I=2}{\cal O}_I\ \amalg\ \coprod_{q\in(\pr{1}\K)'}{\cal O}_q \textrm{ with} \label{eq:24decomp} \\
{\cal O}_I&:=\left\{(\K^2;\K v_1,\K v_2,\K v_3, \K v_4)\in X_g\left|\K v_i=\K v_j \textrm{ for }i,j\in I\right.\right\}, \textrm{ for }I\subset\{1,2,3,4\},\, \#I=2, \notag \\
{\cal O}_q&:=\left\{\left(\K^2;\K v_1,\K v_2,\K (v_1+v_2), \K (q_1v_1+q_2v_2)\right)\in X_g\right\} \textrm{ for }q=\K(q_1e_1+q_2e_2)\in(\pr{1}\K)' \notag 
\end{align}
where $(\pr{1}\K)':=\pr{1}\K\setminus\{\K e_1,\, \K e_2,\, \K(e_1+e_2)\}$ and $\{e_1,e_2\}$ is the standard basis of $\K^2$. 
Now, we can see that a $4$-tuple flag in $X_{2,4}$ is contained in $X_g$ if and only if it is indecomposable by definition of $X_g$ (\ref{eq:xgdef}). 
Furthermore, the indecomposability is equivalent to be contained in the orbits in $\pi^{-1}(p_g)$ by definition of $\pi$ in (\ref{eq:mainsurj}) and the list of elements of $\mathcal{P}=\mathcal{P}_{2,3}$ in (\ref{eq:24ex}). 
Hence, we have $\pi^{-1}(p_g)=\mathrm{diag}(G)\backslash X_g$.
\end{exmp}

\subsection{Representatives of orbits \label{representatives}}

In this section, we give representatives of $\mathrm{diag}(G)$-orbits on $G^m/P^m$, and describe the orbit decomposition of each fibre of the surjection $\pi$ defined in (\ref{eq:mainsurj}). 
For this aim, we give $3$ types of representatives of $m$-tuple flags occuring as summands in the splitting (\ref{eq:splitspecific}). 

\begin{enumerate}
\item

\noindent For the set $\{1,2,\ldots ,k\}\subset\{1,2,\ldots,m\}$ with $0\leq k\leq m$, we define an $m$-tuple flag by 
\begin{subequations}
\begin{equation} \label{eq:repI}
F(\{1,2,\ldots,k\}):=\left(\K;\overbrace{\K,\K,\ldots,\K}^{k},\overbrace{0,0,\ldots,0}^{m-k}\right)
\end{equation}
More generally, for a set $I\subset\{1,2,\cdots ,m\}$, we define an $m$-tuple flag $F(I)$ by
\begin{align}
F(I):=&\left(\K;V^{(1)},V^{(2)},\ldots,V^{(m)}\right), \textrm{ where}  \label{eq:repIp} \\
V^{(i)}=&\begin{cases} \K & \textrm{if }i\in I \\ 0 & \textrm{if }i\notin I \end{cases} \textrm{ for }1\leq i\leq m. \notag
\end{align}
Then we have 
\begin{equation}
F(I)\in Fl_{\bd(I,1)}(\K). \label{eq:repId}
\end{equation}
\end{subequations}

\item

\noindent For the set $\{1,2,\ldots ,r,r+1\}\subset\{1,2,\ldots,m\}$ where $2\leq r\leq m-1$, we define
\begin{subequations}
\begin{equation}
D(\{1,2,\ldots,r,r+1\}):=\left(\K^r;\K e_1,\K e_2,\ldots,\K e_r, \K\left(\sum_{h=1}^re_h\right),0^{m-r-1}\right) \label{eq:repJ}
\end{equation}
where $\{e_i\}_{i=1}^r$ is the standard basis of $\K^{r}$, and $0^l$ denotes the repetition of $0$ for $l$ times. 

More generally, for a set $J=\{j_1< j_2<\cdots <j_r<j_{r+1}\}\subset\{1,2,\ldots ,m\}$ where $2\leq r\leq m-1$, we define an $m$-tuple flag $D(J)$ by
\begin{align}
D(J):=&\left(\K^r;V^{(1)},V^{(2)},\ldots,V^{(m)}\right), \textrm{ where} \label{eq:repJp} \\
V^{(i)}=&\begin{cases} \K e_h & \textrm{if }i=j_h\in J, \;\;\;1\leq h\leq r \\ \K \left(\sum_{h=1}^re_h\right) & \textrm{if }i=j_{r+1}\in J \\ 0 & \textrm{if }i\notin J \end{cases}  \textrm{   for }1\leq i\leq m. \notag
\end{align}
Then we have 
\begin{equation}
D(J)\in Fl_{\bd(J,\#J-1)}(\K^{\#J-1}), \textrm{ and }3\leq \#J. \label{eq:repJd}
\end{equation}
\end{subequations}

\item

\noindent For the set $\{1,2,\ldots ,k\}\subset\{1,2,\ldots,m\}$, an integer $r$ satisfying $4\leq r+2\leq k$, and $q=(q^{(r+2)},q^{(r+3)},\ldots ,q^{(k)})\in(\pr{r-1}\K)^{k-r-1}$, we define
\begin{subequations}
\begin{equation}
E(\{1,\ldots,k\},r,q):=\left(\K^r;\K e_1,\ldots,\K e_r, \K\left(\sum_{h=1}^re_h\right),q^{(r+2)},\ldots,q^{(k)}, 0^{m-k}\right). \label{eq:repK}
\end{equation}

More generally, for a set $K=\{k_1<k_2<\cdots <k_{\#K}\}\subset \{1,2,\ldots ,m\}$, an integer $r$ satisfying $4\leq r+2\leq \#K$, and $q=(q^{(r+2)},q^{(r+3)},\ldots ,q^{(\#K)})\in(\pr{r-1}\K)^{\#K-r-1}$, we define an $m$-tuple flag $E(K,r,q)$ by
\begin{align}
E(K,r,q):=&\left(\K^r;V^{(1)},V^{(2)},\ldots,V^{(m)}\right), \textrm{ where} \label{eq:repKp} \\
V^{(i)}=&\begin{cases} \K e_h & \textrm{if }i=k_h\in K, \;\;\;1\leq h\leq r \\ \K \left(\sum_{h=1}^re_h\right) & \textrm{if }i=k_{r+1}\in K \\ q^{(h)} & \textrm{if }i=k_h\;\;\; r+2\leq h\leq \#K \\ 0 & \textrm{if }i\notin K \end{cases} \textrm{   for }1\leq i\leq m.   \notag
\end{align}
Then  we have
\begin{equation}
E(K,r,q)\in Fl_{\bd(K,r)}(\K^r), \textrm{ and }4\leq r+2\leq \#K. \label{eq:repKd}
\end{equation}
\end{subequations}
\end{enumerate}

All of them are indecomposable from Proposition \ref{thm:indrep}. 
Hence if we define an $m$-tuple flag 
\begin{equation} \label{eq:flagsum}
F(p,q):=F(\emptyset)^{\oplus (n-r(p))} \oplus \left(\bigoplus_{a=1}^A F(I_a)\right) \oplus \left( \bigoplus_{b=1}^B D(J_b)\right) \oplus \left( \bigoplus_{c=1}^C E(K_c,r_c,q_{c}) \right)
\end{equation}
for $p\in\mathcal{P}$ of the form (\ref{eq:pnmelement}) and $q=(q_{1}, q_{2}, \ldots ,q_{C})\in\prod_{c=1}^C(\pr{r_c-1}\K)^{\#K_c-r_c-1}$, 
then from Proposition \ref{thm:pnmid}, $F(p,q)\in Fl_{\bd(\{1,2,\ldots,m\},n)}(\K^n)=G^m/P^m$, and it is $p$-constructible from (\ref{eq:repId}), (\ref{eq:repJd}), and (\ref{eq:repKd}) (see (\ref{eq:splitspecific})).
Hence $\mathrm{diag}(G)\cdot F(p,q)\in\pi^{-1}(p)$ by definition of $\pi$ in (\ref{eq:mainsurj}). 
With these notations, we state the second main result as follows:
\begin{thm} \label{thm:rep} Let $G,\, P,\,\mathcal{P}=\mathcal{P}_{n,m}$, and $\pi$ be as in Theorem \ref{thm:main}. 
For the surjection $\pi: \mathrm{diag}(G)\backslash G^m/P^m\twoheadrightarrow\mathcal{P}$ and an element $p\in\mathcal{P}$ of the form (\ref{eq:pnmelement}), we can define a map as  
\begin{align} \label{eq:emb}
\iota_p\colon\prod_{c=1}^C(\pr{r_c-1}\K)^{\#K_c-r_c-1}\rightarrow\pi^{-1}(p),\, q \mapsto \mathrm{diag}(G)\cdot F(p,q) 
\end{align}
where $F(p,q)\in G^m/P^m$ is the $m$-tuple flag defined in (\ref{eq:flagsum}). 
It is an open dense embedding, and if $C=0$ (in other words, the 3rd family of $p$ is empty), then $\pi^{-1}(p)$ is a singleton and the map (\ref{eq:emb}) gives a bijection between singletons. 
\end{thm}

To explain the meaning of notations, we compare the descriptions of orbits obtained by the direct computations and those obtained by Theorems \ref{thm:main} and \ref{thm:rep} in the following exmaples. 

\begin{exmp}[see Example \ref{thm:23}] \label{thm:23comp}
If $(n,m)=(2,3)$, we have seen that $\pi$ is bijective by listing up all elements of $\mathcal{P}=\mathcal{P}_{2,3}$ and by determining the target of $\pi$ for each orbit obtained by a direct computation. 
Alternatively, we can check that $\pi$ is bijective from Theorem \ref{thm:rep} without the direct computation of the orbit decomposition since the third family of each element of $\mathcal{P}$ is empty in this case. 

Hence each fibre of $p\in\mathcal{P}$ is a singleton and the representative of the corresponding orbit is given as follows by (\ref{eq:emb}):
\begin{align*}
\pi^{-1}(p_c)&=\left\{\mathrm{diag}(G)\cdot F(p_c,\ast)\right\}=\left\{\mathrm{diag}(G)\cdot\left(F(\emptyset)\oplus F(\{1,2,3\})\right)\right\} \\
&=\left\{\mathrm{diag}(G)\cdot(\K^2; \K e_1, \K e_1, \K e_1)\right\} \\
\pi^{-1}(p_1)&=\left\{\mathrm{diag}(G)\cdot F(p_1,\ast)\right\}=\left\{\mathrm{diag}(G)\cdot\left(F(\{1\})\oplus F(\{2,3\})\right)\right\} \\
&=\left\{\mathrm{diag}(G)\cdot(\K^2; \K e_1, \K e_2, \K e_2)\right\} \\
\pi^{-1}(p_o)&=\left\{\mathrm{diag}(G)\cdot F(p_o,\ast)\right\}=\left\{\mathrm{diag}(G)\cdot D(\{1,2,3\})\right\} \\
&=\left\{\mathrm{diag}(G)\cdot(\K^2; \K e_1, \K e_2, \K (e_1+e_2))\right\} 
\end{align*}
where $\{e_1,e_2\}$ is the standard bases in $\K^2$. 
Similarly to the case for $p_1$, we have $\pi^{-1}(p_2)=\{\mathrm{diag}(G)\cdot(\K^2; \K e_2, \K e_1, \K e_2)\}$, $\pi^{-1}(p_3)=\{\mathrm{diag}(G)\cdot(\K^2; \K e_2, \K e_2, \K e_1)\}$.
\end{exmp}

\begin{exmp}[see Examples \ref{thm:hausdorff} and \ref{thm:24}] \label{thm:24comp}
Let $(n,m)=(2,4)$. 
For $p_g\in\mathcal{P}=\mathcal{P}_{2,4}$ defined in (\ref{eq:24ex}), the open dense embedding $\iota_{p_g}\colon \pr{1}\K\hookrightarrow\pi^{-1}(p_g)$ is given by
\begin{equation*}
q\mapsto \mathrm{diag}(G)\cdot F(p_g,q)=\mathrm{diag}(G)\cdot(\K^2; \K e_1, \K e_2, \K (e_1+e_2), q)=\begin{cases} {\cal O}_q & \textrm{if } q\in(\pr{1}\K)' \\ {\cal O}_{\{1,4\}} & \textrm{if } q=\K e_1 \\ {\cal O}_{\{2,4\}} & \textrm{if } q=\K e_2 \\ {\cal O}_{\{3,4\}} & \textrm{if } q=\K (e_1+e_2).\end{cases}
\end{equation*}
The complement of the image of this embedding into $\pi^{-1}(p_g)$ is a three-point set
\begin{equation*}
\left\{{\cal O}_{\{1,2\}}, {\cal O}_{\{2,3\}}, {\cal O}_{\{3,1\}}\right\}.
\end{equation*} 
\end{exmp}

\section{Proof of the main results}

To prove Theorems \ref{thm:main} and \ref{thm:rep}, we formulate the notion of splitting of multiple flags into indecomposable ones in Section 3.1, and introduce some remarkable previous results on this notion according to \cite{mwzA} in Section 3.2. 
With these preparations, we divide the main results into $3$ essential parts (Propositions \ref{thm:inddimclasscomp}, \ref{thm:indrep}, and \ref{thm:pnmid}), and give proofs of them in Sections 3.3 and 3.4. 

\subsection{Splitting of multiple flags into indecomposables}

In this section, we introduce the following notations to formulate the notion of splitting of multiple flags into indecomposable ones. 
\begin{df} \label{thm:dimvec}
Let $p_1$, $p_2,\ldots,\,p_m$ be positive integers, then we define a set $\Delta_{p_1,p_2,\ldots,p_m}$ by
\begin{align*}
\Delta_{p_1,p_2,\ldots,p_m}:=\left\{\left(\ba^{(1)}, \ba^{(2)}, \cdots , \ba^{(m)}\right)\in \N^{p_1}\times\N^{p_2}\times\cdots\times\N^{p_m} \right. & \\
\left. \left| |\ba^{(1)}|=|\ba^{(2)}|=\cdots =|\ba^{(m)}|\geq 1 \right. \right\}&
\end{align*}
where \[|\ba^{(i)}|=\sum_{j=1}^{p_i}a_j^{(i)} \textrm{  for  } \ba^{(i)}=\fourIdx{t}{}{}{}{\left(a_1^{(i)}, a_2^{(i)}, \cdots , a_{p_i}^{(i)} \right)\in\N^{p_i}},\;\;1\leq i\leq m.\]
For an element $\ba=\left(\ba^{(1)}, \ba^{(2)}, \cdots , \ba^{(m)}\right)\in\Delta_{p_1,p_2,\ldots,p_m}$, we define
\[|\ba|:= |\ba^{(1)}|=|\ba^{(2)}|=\cdots =|\ba^{(m)}|.\]
\end{df}
In other words, $\Delta_{p_1,p_2,\ldots,p_m}$ is the set of all $m$-tuples $\ba$ of compositions of the common positive integers $|\ba|$ with the lengths $p_1$, $p_2,\ldots,\,p_m$. 
From now on, we call the common number $|\ba|$ the size of $\ba$, and call an element of $\Delta_{p_1,p_2,\ldots,p_m}$ an abstract dimension vector. 

\bigskip
Next, we introduce an additive category ${\cal F}_{p_1,p_2,\ldots,p_m}$. 

\noindent \underline{Objects}: 
An object of ${\cal F}_{p_1,p_2,\ldots,p_m}$ is a tuple $(V; F^{(1)}, F^{(2)}, \ldots , F^{(m)})$ where $V$ is a vector space and each $F^{(i)}$ is a sequence of subspaces of $V$ such as
\[
F^{(i)}=\bigl(  0=V_0^{(i)}\subset V_1^{(i)}\subset \cdots \subset V_{p_i-1}^{(i)}\subset V_{p_i}^{(i)}=V \bigr) \textrm{  for  }1\leq i\leq m.
\]
We call an object $(V; F^{(1)}, F^{(2)}, \ldots , F^{(m)})$ of ${\cal F}_{p_1,p_2,\ldots,p_m}$ an $m$-tuple flag, and $V$ its whole vector space. 

\noindent \underline{Morphisms}: 

For $m$-tuple flags $F$ and $G$ which are objects of ${\cal F}_{p_1,p_2,\ldots,p_m}$ of the forms
\begin{gather}
F=(V;F^{(1)}, F^{(2)}, \ldots , F^{(m)}),\;\;G=(W;G^{(1)}, G^{(2)}, \ldots , G^{(m)}) \textrm{ where}  \label{eq:flagcat} \\
\begin{array}{rcl}
F^{(i)}&=&\bigl(  0=V_0^{(i)}\subset V_1^{(i)}\subset \cdots \subset V_{p_i-1}^{(i)}\subset V_{p_i}^{(i)}=V \bigr), \\
G^{(i)}&=&\bigl(  0=W_0^{(i)}\subset W_1^{(i)}\subset \cdots \subset W_{p_i-1}^{(i)}\subset W_{p_i}^{(i)}=W \bigr)
\end{array}
\textrm{  for  }1\leq i\leq m, \notag
\end{gather}
a morphism $f$ from $F$ to $G$ is a linear map from $V$ to $W$ with the condition 
\[
f(V_{j}^{(i)})\subset W_j^{(i)} \textrm{  for }1\leq i\leq m,\; 1\leq j\leq p_i. 
\]


\noindent \underline{Direct sums}: 
We define the notion of the direct sum of $m$-tuple flags as the collection of direct sums of each subspaces. 
More precisely, given $m$-tuple flags $F$ and $G$ which are objects of ${\cal F}_{p_1,p_2,\ldots,p_m}$ as in (\ref{eq:flagcat}), the direct sum $F\oplus G$ which is also an obeject of ${\cal F}_{p_1,p_2,\ldots,p_m}$ is defined as 
\begin{align}
F\oplus G&:=(V\oplus W; F^{(1)}\oplus G^{(1)}, F^{(2)}\oplus G^{(2)}, \ldots , F^{(m)}\oplus G^{(m)}) \textrm{  where} \label{eq:directsum} \\ 
F^{(i)}\oplus G^{(i)}&:=\bigl( 0=V_0^{(i)}\oplus W_0^{(i)} \subset V_1^{(i)}\oplus W_1^{(i)}\subset \cdots \subset V_{p_i}^{(i)}\oplus W_{p_i}^{(i)}=V\oplus W \bigr)\textrm{  for  }1\leq i\leq m. \notag
\end{align} 
With these definitions, ${\cal F}_{p_1,p_2,\ldots,p_m}$ is an additive category. 
We say that an $m$-tuple flag (an object of ${\cal F}_{p_1,p_2,\ldots,p_m}$) is \emph{indecomposable} if it cannot be realised as the direct sum of two non-trivial (the whole space is non-zero) $m$-tuple flags. 
%
For an $m$-tuple flag $F$ which is an object of ${\cal F}_{p_1,p_2,\ldots,p_m}$ of the form (\ref{eq:flagcat}), we define $\dim F\in\Delta_{p_1,p_2,\ldots,p_m}$ by 
\begin{align*}
\dim F :=&\left(\ba^{(1)},\ba^{(2)},\ldots,\ba^{(m)}\right)\in\N^{p_1}\times\N^{p_2}\times\cdots\times\N^{p_m} \textrm{  where} \\
\ba^{(i)}=&\left(\dim V_j^{(i)}/V_{j-1}^{(i)}\right)_{j=1}^{p_i}\in\N^{p_i} \textrm{  for }1\leq i\leq m.
\end{align*}
We call it the dimension vector of $F$. 
Remark that for isomorphism classes $\lambda=\left[F_{\lambda}\right]$ and $\mu=\left[F_{\mu}\right]$ of objects of ${\cal F}_{p_1,\ldots,p_m}$, the dimension vector $\dim F_{\lambda}\in\Delta_{p_1,p_2,\ldots,p_m}$ and the isomorphism class $\left[F_{\lambda}\oplus F_{\mu}\right]$ do not depend on the choice of representatives of isomorphism classes. 
Hence we write them just $\dim\lambda\in\Delta_{p_1,p_2,\ldots,p_m}$ and $\lambda\oplus\mu$. 
Then we have 
\begin{gather}
\dim(\lambda\oplus\mu)=\dim\lambda+\dim\mu. \label{eq:dimdirectsum} 
\end{gather}

From now on, we fix $p_1,p_2,\ldots,p_m$ and write just $\Delta$ and ${\cal F}$ for $\Delta_{p_1,p_2,\ldots,p_m}$ and ${\cal F}_{p_1,p_2\ldots,p_m}$. 

Let us fix an abstract dimension vector $\ba=(\ba^{(i)})_{i=1}^m\in\Delta$ and an $|\ba|$-dimensional vector space $V$. 
We define $Fl_{\ba}(V)$ as the set of all objects $F$ of ${\cal F}$ whose whole vector spaces are $V$ and $\dim F=\ba$. 
Let $P_{\ba^{(i)}}$ be the parabolic subgroup of $GL(V)$ with the Jordan blocks of the sizes $a_1^{(i)}, a_2^{(i)}, \ldots, a_{p_i}^{(i)}$ for $1\leq i\leq m$ where $\ba^{(i)}=(a_j^{(i)})_{j=1}^{p_i}$. 
We set $P_{\ba}:=\prod_{i=1}^mP_{\ba^{(i)}}$.  
Then the $m$-tuple flag variety $GL(V)^m/P_{\ba}$ which we are dealing with in this article is naturally identified with $Fl_{\ba}(V)$, and the identification is $\mathrm{diag}(GL(V))$-equivariant. 
Furthermore, it is equivalent for two $m$-tuple flags with the same whole vector space $V$, that they are contained in the same $\mathrm{diag}(GL(V))$-orbit and that they are isomorphic as objects of ${\cal F}$. 
Hence we have a natural identification
\begin{equation} \label{eq:Flcat}
\mathrm{diag}(GL(V))\backslash GL(V)^m/P_{\ba}\cong \left\{F\in\mathrm{ob}({\cal F})\;|\; \dim F=\ba \right\}/\underset{\mathrm{isom}}{\sim}.
\end{equation}


To describe isomorphism classes of $m$-tuple flags, we introduce the followings: 
\begin{align}
\tilde{\Lambda}_{\ba}&:=\left\{I\in\mathrm{ob}({\cal F})\left|\dim I\leq \ba,\textrm{ and } I \textrm{ is indecomposable} \right.\right\}/\underset{\mathrm{isom}}{\sim}; \label{eq:tildelambda} \\
\Lambda_{\ba}&:= \left\{\dim\lambda\in\Delta \left| \lambda\in\tilde{\Lambda}_{\ba}\right.\right\};  \label{eq:lambda} \\
Fl^{\mathrm{ind}}_{\bd}(W)&:= \left\{I\in Fl_{\bd}(W)\left| [I]\in\tilde{\Lambda}_{\ba} \right.\right\} \textrm{ for }\bd\in\Lambda_{\ba}\textrm{ and a }|\bd|\textrm{-dimensional vector spcae }W. \label{eq:Flind}
\end{align}
A remarkable property of ${\cal F}$ is that it is a Krull-Schmidt category \cite{mwzA}: any object of ${\cal F}$ always has a unique splitting into indecomposable ones up to isomorphisms. 
Hence, we have the following natural bijection:
\begin{equation} \label{eq:split}
\begin{array}{ccc} \left\{F\in\mathrm{ob}({\cal F})\;|\; \dim F=\ba \right\}/\underset{\mathrm{isom}}{\sim} & \xrightarrow{\sim} & \displaystyle{\left\{ \left(m_{\lambda}\right)_{\lambda\in\tilde{\Lambda}_{\ba}}\in \N^{\tilde{\Lambda}_{\ba}} \left| \sum_{\lambda\in\tilde{\Lambda}_{\ba}}m_{\lambda}\dim \lambda= \ba \right.\right\}} \\
\rotatebox{90}{$\in$} & \null & \rotatebox{90}{$\in$} \\ \vspace{2mm}
\bigoplus_{\lambda\in\tilde{\Lambda}_{\ba}}\lambda^{\oplus m_{\lambda}} & \mapsto & \left(m_{\lambda}\right)_{\lambda\in\tilde{\Lambda}_{\ba}}. \end{array}
\end{equation}
Furthermore, we have a natural surjection
\begin{equation} \label{eq:surjection}
\begin{array}{ccc}  \displaystyle{\left\{ \left(m_{\lambda}\right)_{\lambda\in\tilde{\Lambda}_{\ba}}\in\N^{\tilde{\Lambda}_{\ba}} \left| \sum_{\lambda\in\tilde{\Lambda}_{\ba}}m_{\lambda}\dim\lambda= \ba \right.\right\}} & \twoheadrightarrow & \displaystyle{\left\{ \left(m_{\bd} \right)_{\bd\in\Lambda_{\ba}}\in\N^{\Lambda_{\ba}} \left| \sum_{\bd\in\Lambda_{\ba}}m_{\bd}\bd= \ba \right.\right\}} \\
\rotatebox{90}{$\in$} & \null & \rotatebox{90}{$\in$} \\ \vspace{2mm}
\left(m_{\lambda}\right)_{\lambda\in\tilde{\Lambda}_{\ba}} & \mapsto & \left( \sum_{\lambda\in\tilde{\Lambda}_{\ba},\,\dim\lambda=\bd} m_{\lambda} \right)_{\bd \in \Lambda_{\ba} } \end{array}
\end{equation}
by forgetting differences among distinct isomorphism classes of indecomposable $m$-tuple flags with the same dimension vectors. 
Let $\mathcal{M}_{\ba}$ denote the set in the right-hand side of (\ref{eq:surjection}). 
Remark that $\mathcal{M}_{\ba}$ is a finite set by definition. 
Concequently, we have the following surjection as the composition of the bijections (\ref{eq:Flcat}) and (\ref{eq:split}) and the surjection (\ref{eq:surjection}): 
\begin{gather}
\begin{array}{ccc}  \mathrm{diag}(GL(V))\backslash GL(V)^m/P_{\ba} & \twoheadrightarrow & \mathcal{M}_{\ba} \\
\rotatebox{90}{$\in$} & \null & \rotatebox{90}{$\in$} \\ \vspace{2mm}
\mathrm{diag}(GL(V))\cdot F & \mapsto & \left(m_{\bd}\right)_{\bd\in\Lambda_{\ba}} \end{array} \textrm{ where} \label{eq:surjabstract} \\
F\underset{\mathrm{isom}}{\sim}\bigoplus_{\bd\in\Lambda_{\ba}}\bigoplus_{i=1}^{m_{\bd}}F_{\bd,i}, \textrm{ and }
F_{\bd,i}\in Fl_{\bd}^{\mathrm{ind}}(\K^{|\bd|}) \textrm{ for }\bd\in\Lambda_{\ba},\,1\leq i\leq m_{\bd}. \notag
\end{gather}
The fibre of the surjection (\ref{eq:surjabstract}) at $(m_{\bd})_{\bd\in\Lambda_{\ba}}\in\mathcal{M}_{\ba}$ is homeomorphic to
\begin{equation} \label{eq:preimage}
\prod_{\bd\in\Lambda_{\ba}}\left(\mathrm{diag}(\GL{|\bd|})\backslash Fl_{\bd}^{\mathrm{ind}}(\K^{|\bd|})\right)^{m_{\bd}}/\mathcal{S}_{m_{\bd}}. 
\end{equation}
Remark that if $\mathrm{diag}(\GL{|\bd|})\backslash Fl_{\bd}^{\mathrm{ind}}(\K^{|\bd|})$ is a singleton (indecomposable $m$-tuple flag with the dimension vector $\bd$ is unique up to isomorphisms) for all $\bd$ with $m_{\bd}\geq 1$, then the fibre is also a singleton. 

Considering the surjection (\ref{eq:surjabstract}), we can separate the proof of the description of $\mathrm{diag}(GL(V))$-orbits on $Fl_{\ba}(V)\cong GL(V)^m/P_{\ba}$ into the following 3 steps. 
\begin{itemize}
\item[Step1] Describing the set $\Lambda_{\ba}$ defined in (\ref{eq:lambda}): classifying all dimension vectors $\bd\leq\ba$ which has at least one indecomposable $m$-tuple flag.
\item[Step2] Describing the set $\mathcal{M}_{\ba}$ which is the target space of the surjection (\ref{eq:surjabstract}): determining all combinations of dimension vectors in $\Lambda_{\ba}$ whose summation coincides with $\ba$, which is a purely combinatorial problem. 
\item[Step3] Describing the $\mathrm{diag}(\GL{|\bd|})$-orbit decomposition of $Fl_{\bd}^{\mathrm{ind}}(\K^{|\bd|})$ introduced in (\ref{eq:Flind}) for each $\bd\in\Lambda_{\ba}$ which parametrises fibres of the surjection (\ref{eq:surjabstract}). 
\end{itemize}

\subsection{Tits forms, finiteness and indecomposability}

There are some previous results on these steps which we focus on in this section. 
To classify dimension vectors, we introduce some notions for them. 
\begin{df}
Let $\bd\in\Delta=\Delta_{p_1,p_2,\ldots,p_m}$ (see Definition \ref{thm:dimvec}), then 
\begin{enumerate}
\item $\bd$ is said to be finite (resp. infinite) if $\mathrm{diag}(\GL{|\bd|})\backslash Fl_{\bd}(\K^{|\bd|})$ (the set of all isomorphism classes of objects of ${\cal F}={\cal F}_{p_1,p_2,\ldots,p_m}$ with the dimension vector $\bd$) is finite (resp. infinite); 
\item $\bd$ is said to be indecomposable (resp. decomposable) if there exist some (resp. no) indecomposable objects of ${\cal F}$ with the dimension vector $\bd$. 
\end{enumerate}
\end{df}
Under this notation, we can define $\Lambda_{\ba}$ alternatively as $\{\bd\in\Delta\left|\bd\leq\ba\textrm{ and indecomposable}\right.\}$. 
To determine an abstract dimension vector indecomposable or not, we define a quadratic form called Tits form. 
\begin{df}[{\cite[(2.1)]{mwzA}}] \label{thm:tits}
For an abstract dimension vector $\bd\in\Delta$, we define
\[
Q(\bd):= \dim GL(V)-\dim FL_{\bd}(V)
\]
where $V$ is a $|\bd|$-dimensional vector space, and call it the Tits form of $\bd$.  
\end{df}
\noindent 
An easy computation leads that 
\[
Q(\bd)=\frac{2-m}{2}|\bd|^2+\frac{\|\bd\|^2}2
\]
where $\|\bd\|^2=\sum_{i=1}^m{}^t\bd^{(i)}\bd^{(i)}$ for $\bd=(\bd^{(1)},\bd^{(2)},\ldots,\bd^{(m)})$. 


By introducing these notations, relationships between finiteness, indecomposability, and Tits form are given as follows: 
\begin{fact}[{\cite[Prop. 3.3]{mwzA}}] \label{thm:findim}
If an abstract dimension vector $\bd$ is finite, then there is an open $\mathrm{diag}(GL(V))$-orbit on $Fl_{\bd}(V)$. 
Furthermore if $Fl_{\bd}(V)$ has an open $\mathrm{diag}(GL(V))$-orbit, then $Q(\bd)\geq 1$. 
In particular, if an abstract dimension vector $\bd$ is finite, then $Q(\bd)\geq 1$. 
\end{fact}
\begin{fact}[{\cite[Thm. 1]{kac}}, {\cite[Prop. 3.1]{mwzA}}] \label{thm:inddim}
Let $\bd$ be an abstract dimension vector. 
If $\bd$ is indecomposable, then $Q(\bd)\leq 1$. 
Furthermore, if $\bd$ is indecomposable and $Q(\bd)=1$, then indecomposable multiple flags with the dimension vector $\bd$ are unique up to isomorphisms. 
\end{fact}
\noindent In other words, the consequence of the second claim says that $\#\mathrm{diag}(\GL{|\bd|})\backslash Fl_{\bd}^{\mathrm{ind}}(\K^{|\bd|})=1$ for an indecomposable dimension vector $\bd$ such that $Q(\bd)=1$. 


\begin{remark}
From these results, \cite{mwzA} proved that the surjection (\ref{eq:surjabstract}) is in fact a bijection if $\ba$ is finite. 
On the other hand, we are dealing with infinite dimension vectors in this article, hence the surjection (\ref{eq:surjabstract}) is no more injective. 
\end{remark}

\subsection{Indecomposable dimension vectors and representatives of multiple flags}

Let us go back to our setting. 
We fix $G=\GL{n}$, $\K=\C$ or $\R$, and $P$ denotes the maximal parabolic subgroup of $G$ such that $G/P\cong\pr{n-1}\K$.
Our main problem is to describe $\mathrm{diag}(G)$-orbit decomposition of the multiple flag variety $G^m/P^m$. 
Let $\ba_{n,m}:=({}^t(1,n-1)^m)\in\Delta_{2^m}$ where ${}^t(1,n-1)^m$ (resp. $2^m$) denotes the repretition of ${}^t(1,n-1)$ (resp. $2$) for $m$ times. 
From now on, we fix $n$, $m$, and $p_1=p_2=\cdots=p_m=2$, and we write $\Delta$, ${\cal F}$, and $\ba$ instead of $\Delta_{2^m}$, ${\cal F}_{2^m}$, and $\ba_{n,m}$. 
Furthermore, we write $(|\bd|; d^{(1)}_1, d^{(2)}_1, \ldots, d^{(m)}_1)$ instead of an abstract dimension vector $\bd=(\bd^{(1)},\bd^{(2)},\ldots,\bd^{(m)})\in\Delta=\Delta_{2^m}$ where $\bd^{(i)}={}^t(d^{(i)}_1,d^{(i)}_2)={}^t(d^{(i)}_1,|\bd|-d^{(i)}_1)$ for $1\leq i\leq m$ for simplicity. 
Under this notation, summation defined in Section 2.1 is compatible with that of $\Delta$. 
Furthermore, a tuple of integers $\bd(K,r)=(r;\delta_{1,K},\delta_{2,K},\ldots,\delta_{m,K})$ defined in (\ref{eq:dkr}) for some $K\subset\{1,2,\ldots,m\}$ and $r\in\Z_{+}$ is an element of $\Delta$, and $\ba=\ba_{n,m}$ is written as $(n;1^m)=\bd(\{1,2,\ldots,m\},n)$ defined in (\ref{eq:dkr}). 
According to the notations in Section 3.1, $P^m=P_{\ba_{n,m}}\subset G^m$, and we have $G^m/P^m\cong Fl_{\ba_{n,m}}(\K^n)$. 


If an abstract dimension vector $\bd\in\Delta$ is bounded by $\ba=\ba_{n,m}=\bd(\{1,2,\ldots,m\},n)$, then it is of the form $\bd(K, r)$ defined in (\ref{eq:dkr}) for some $K\subset\{1,2,\ldots,m\}$ and $r\in\N$. 
By Definition \ref{thm:tits}, we can compute that
\[
Q(\bd(K,r))=(r-1)(r-\#K+1)+1.
\]
From Fact \ref{thm:inddim} we have the following proposition. 
\begin{prop} \label{thm:inddimclass}
For abstract dimension vectors $\bd(K,r)\in\Delta=\Delta_{2^m}$ defined in (\ref{eq:dkr}), we have the following properties: 
\begin{enumerate}
\item For the set $\Lambda_{\ba}$ (see (\ref{eq:lambda})) with $\ba=\ba_{n,m}\in\Delta$, we have
\begin{align}
\Lambda_{\ba}\subset&\left\{\bd(K,r)\leq\ba\left|K\subset\{1,2,\ldots,m\},\,r\in\N,\,\textrm{ and }r=1\textrm{ or }3\leq r+1\leq\#K\right.\right\} \label{eq:Lambda} \\
=& \left\{\bd(I,1)\leq\ba\left|I\subset\{1,2,\ldots,m\}\right.\right\} \notag \\
&\;\amalg\left\{\bd(J,\#J-1)\leq\ba\left|J\subset\{1,2,\ldots,m\},\,3\leq \#J\right.\right\} \notag \\
&\;\amalg\left\{\bd(K,r)\leq\ba\left|K\subset\{1,2,\ldots,m\},\,r\in\N,\,4\leq r+2\leq \#K\right.\right\}.  \notag
\end{align}
\item The set $\mathrm{diag}(\GL{|\bd|})\backslash Fl_{\bd}^{\mathrm{ind}}(\K^{|\bd|})$ (see (\ref{eq:Flind})) is a singleton or empty if $\bd=\bd(I,1)$ or $\bd(J,\#J-1)\in\Delta$ for some $I, J\subset\{1,2,\ldots,m\}$ satisfying $3\leq \#J$. 
\end{enumerate}
\end{prop}
Motivated by this proposition, to complete the determination of the set $\Lambda_{\ba}$ (the set of all indecomposable dimension vectors $\bd\in\Delta$ satisfying $\bd\leq\ba$: the first step of our proof we introduced in Section 3.1) and to describe the orbit space $\mathrm{diag}(\GL{|\bd|})\backslash Fl_{\bd}^{\mathrm{ind}}(\K^{|\bd|})$ (the set of all isomorphism classes of indecomposable $m$-tuple flags whose dimension vectors are $\bd$: the third step of our proof we introduced in Section 3.1) for $\bd\in\Lambda_{\ba}$, we introduce some indecomposable $m$-tuple flags explicitly for dimension vectors $\bd(K,r)$ in the right-hand side of (\ref{eq:Lambda}). 
\begin{enumerate}
\item 
Let $I$ be a subset of $\{1,2,\ldots,m\}$. 
Proposition \ref{thm:inddimclass} says that there is at most one isomorphism class of indecomposable $m$-tuple flags having the dimension vector $\bd(I,1)$. 
In fact, the $m$-tuple flag $F(I)$ defined in (\ref{eq:repIp}) is indecomposable since $1$-dimensional vector space does not have any decompositions into non-zero subspaces. 
As we saw in (\ref{eq:repId}), $\dim F(I)=\bd(I,1)$. 
\item 
Let $J=\{j_1<j_2<\cdots<j_r<j_{r+1}\}\subset\{1,2,\ldots,m\}$ satisfy $3\leq\#J=r+1$ and $\bd(J,\#J-1)=\bd(J,r)\leq\ba$. 
Proposition \ref{thm:inddimclass} says that there is at most one isomorphism class of indecomposable $m$-tuple flags having the dimension vector $\bd(J,r)$. 
Consider the $m$-tuple flag $D(J)\in Fl_{\bd(J,r)}(\K^r)$ defined in (\ref{eq:repJp}). 
If $D(J)=(\K^r;V^{(1)},\ldots,V^{(m)})$ is decomposable, then there exist non-zero subspaces $V_1$ and $V_2$ of $\K^r$ such that $\K^r=V_1\oplus V_2$, and $V^{(i)}$ is contained in either $V_1$ or $V_2$ for each $1\leq i\leq m$. 
Then we have a partition $\{1,2,\ldots,r\}=I_1\amalg I_2$ such that $h\in I_k$ if and only if $\K e_h=V^{(j_h)}\subset V_k$ for $1\leq h\leq r$ and $k=1,2$. 
If $I_1$ (resp. $I_2$) is empty, then $V_2$ (resp. $V_1$) coincides with $\K^r$. 
Hence we can assume that $I_1,\,I_2\neq\emptyset$, which leads that $\K\sum_{h=1}^re_h=V^{(j_{r+1})}$ is contained in neither $V_1$ nor $V_2$. 
It contradicts the assumption, hence $D(J)$ is indecomposable. 
As we saw in (\ref{eq:repJd}), $\dim D(J)=\bd(J,r)=\bd(J,\#J-1)$. 
\item 
Let $K=\{k_1<k_2<\cdots <k_{\#K}\}\subset\{1,2,\ldots,m\}$ and $r\in\N$ satisfy $4\leq r+2\leq \#K$ and $\bd(K,r)\leq\ba$. 
For an $m$-tuple flag $(\K^r;V^{(1)},V^{(2)},\ldots,V^{(m)})\in Fl_{\bd(K,r)}(\K^r)$, remark that $V^{(k)}$ is a $\delta_{i,K}$-dimensional subspace in $\K^r$ for each $1\leq i\leq m$. 
Now we define an open dense subset $Fl_{\bd(K,r)}^{\mathrm{gen}}(\K^r)$ of $Fl_{\bd(K,r)}(\K^r)$ as the set of all such $m$-tuple flags such that $\{V^{(k_h)}\}_{h=1}^{r+1}$ is in general position in $\K^r$. 
Then it is contained in $Fl_{\bd(K,r)}^{\mathrm{ind}}(\K^r)$ by a similar arguement as the proof for the indecomposability of $D(J)$. 
Hence $Fl_{\bd(K,r)}^{\mathrm{gen}}(\K^r)$ is an open dense subset of $Fl_{\bd(K,r)}^{\mathrm{ind}}(\K^r)$. 

Now, we have an contineous map $g$ by
\begin{equation*}
g\colon Fl_{\bd(K,r)}^{\mathrm{gen}}(\K^r)\rightarrow PGL(r,\K)
\end{equation*}
such that for $F=(\K^r, V^{(i)})_{i=1}^m\in Fl_{\bd(K,r)}^{\mathrm{gen}}(\K^r)$, $g(F)$ is the unique element of $PGL(r,\K)=\GL{r}/\K^{\times}{}^{\curvearrowright}\pr{r-1}\K$ which transforms $\K e_h$ to $V^{(k_h)}$ for each $1\leq h\leq r$ and $\K( \sum_{h=1}^re_h)$ to $V^{(k_{r+1})}$ where $\{e_i\}_{i=1}^r$ is the standard basis of $\K^r$. 
Since $g(\mathrm{diag}(h)\cdot F)=[h]g(F)\in PGL(r,\K)$ for $h\in\GL{r}$, we can define a contineous map 
\begin{equation} \label{eq:embspecificinv}
\kappa_{\bd(K,r)}\colon \mathrm{diag}(\GL{r})\backslash Fl_{\bd(K,r)}^{\mathrm{gen}}(\K^r) \rightarrow (\pr{r-1}\K)^{\#K-r-1}
\end{equation}
mapping $\mathrm{diag}(\GL{r})\cdot F$ to $(g(F)^{-1}V^{(k_{h})})_{h=r+2}^{\#K}$. 
The inverse of the map (\ref{eq:embspecificinv}) is given by a contineous map $(\pr{r-1}\K)^{\#K-r-1}\ni q\mapsto\mathrm{diag}(\GL{r})\cdot E(K,r,q)$ where $E(K,r,q)$ is the $m$-tuple flag defined in (\ref{eq:repKp}). 
Hence (\ref{eq:embspecificinv}) is a homeomorphism. 
\end{enumerate}
From these arguments, we obtained the results for the first and third step of our proof we introduced in Section 3.1. 
\begin{prop} \label{thm:inddimclasscomp}
The equality holds for (\ref{eq:Lambda}). 
\end{prop}
\begin{prop} \label{thm:indrep} Let $I,\,J,\,K\subset\{1,2,\ldots,m\}$, and $r\in\N$ satisfy $\bd(I,1),\, \bd(J,\#J-1),\, \bd(K,r)\leq\ba=\ba_{n,m}$, $3\leq \#J$, and $4\leq r+2\leq \#K$. 
\begin{enumerate}
\item Indecomposable $m$-tuple flags with the dimension vectors $\bd(I,1)$ and $\bd(J,\#J-1)$ are unique up to isomorphisms, and their representatives are given by $F(I)$ and $D(J)$ respectively which are defined in (\ref{eq:repIp}) and (\ref{eq:repJp}).  
\item For the orbit space $\mathrm{diag}(\GL{r})\backslash Fl_{\bd(K,r)}^{\mathrm{ind}}(\K^r)$ which is naturally identified with the set of all isomorphism classes of indecomposable $m$-tuple flags with the dimension vector $\bd(K,r)$ via the bijection (\ref{eq:Flcat}), we have the following open dense embedding: 
\begin{equation} \label{eq:embspecific}
\iota_{\bd(K,r)}\colon(\pr{r-1}\K)^{\#K-r-1}\hookrightarrow \mathrm{diag}(\GL{r})\backslash Fl_{\bd(K,r)}^{\mathrm{ind}}(\K^r),\,\, q \mapsto \mathrm{diag}(\GL{r})\cdot E(K,r,q)
\end{equation}
where each $E(K,r,q)$ is the $m$-tuple flag defined in (\ref{eq:repKp}). 
\end{enumerate}
\end{prop}

\subsection{Multiplicity of each indecomposable multiple flag}

The remaining part of the proof of Theorems \ref{thm:main} and \ref{thm:rep} is to determine the finite set $\mathcal{M}_{\ba}$ (the set of all multiplicity matrices $\left(m_{\bd}\right)_{\bd\in\Lambda_{\ba}}$ satisfying $\sum m_{\bd}\bd=\ba=\ba_{n,m}$: the second step we introduced in Section 3.1). 
All elements of $\Lambda_{\ba}$ are determined in Proposition \ref{thm:inddimclasscomp}. 
Now, let us consider a linear combination of elements of $\Lambda_{\ba}$ as
\begin{equation} \label{eq:dimsplit}
m_{\bd(\emptyset,1)}\bd(\emptyset,1)+\sum_{a=1}^Am_{\bd(I_a,1)}\bd(I_a,1)+\sum_{b=1}^Bm_{\bd(J_b,\#J_b-1)}\bd(J_b,\#J_b-1)+\sum_{c=1}^Cm_{\bd(K_c,r_c)}\bd(K_c,r_c)
\end{equation}
where $I_a,\, J_b,\, K_c\subset\{1,2,\ldots,m\}$, and $r_c\in\N$ satisfy $m_{\bd(\emptyset,1)},\,m_{\bd(I_a,1)},\,m_{\bd(J_b,\#J_b-1)},\,m_{\bd(K_c,r_c)}\neq0$, $1\leq \#I_a$, $3\leq \#J_b$, and $4\leq r_c+2\leq \#K_c$ for $1\leq a\leq A,\,1\leq b\leq B,\, 1\leq c\leq C$, and each of $\{I_a\}_{a=1}^A,\, \{J_b\}_{b=1}^B$ and $\{K_c\}_{c=1}^C$ is distinct. 
Then (\ref{eq:dimsplit}) is of the form $(r;d^{(1)},\ldots,d^{(m)})$ where
\begin{align}
r&=m_{\bd(\emptyset,1)}+\sum_{a=1}^Am_{\bd(I_a,1)}+\sum_{b=1}^Bm_{\bd(J_b,\#J_b-1)}(\#J_b-1)+\sum_{c=1}^Cm_{\bd(K_c,r_c)}r_c \label{eq:dimsplit1} \\
d^{(i)}&=\sum_{a=1}^Am_{\bd(I_a,1)}\delta_{i,I_a}+\sum_{b=1}^Bm_{\bd(J_b,\#J_b-1)}\delta_{i,J_b}+\sum_{c=1}^Cm_{\bd(K_c,r_c)}\delta_{i,K_c} \textrm{  for }1\leq i\leq m \label{eq:dimsplit2}
\end{align}
from (\ref{eq:dkr}). 
The number (\ref{eq:dimsplit2}) coincides with $1$ if and only if all multiplicities are $1$ and $\{I_a,J_b,K_c\}_{a=1,b=1,c=1}^{A,B,C}$ is a partition of the set $\{1,2,\ldots,m\}$. 
Then the number (\ref{eq:dimsplit1}) coincides with $m_{\bd(\emptyset,1)}+r(p)$ where $r(p)$ is the number introduced in Definition \ref{thm:pnm}. 
Hence, the linear combination (\ref{eq:dimsplit}) coinsides with $\ba=(n;1^m)$ if and only if \[p:=\left(\left\{I_a\right\}_{a=1}^{A}, \left\{J_b\right\}_{b=1}^{B}, \left\{ (K_c,r_c) \right\}_{c=1}^{C} \right)\] is an element of $\mathcal{P}$, $m_{\bd(\emptyset,1)}=n-r(p)$, and $m_{\bd(I_a,1)}=m_{\bd(J_b,\#J_b-1)}=m_{\bd(K_c,r_c)}=1$ for $1\leq a\leq A,\, 1\leq b\leq B$, and $1\leq c\leq C$. 
Hence, we obtain the result for the second step in Section 3.1. 
\begin{prop} \label{thm:pnmid}
For the finite sets $\mathcal{P}=\mathcal{P}_{n,m}$ (see Definition \ref{thm:pnm}) and $\mathcal{M}_{\ba}$ which is the set of all splittings of $\ba=\ba_{n,m}$ into elements of $\Lambda_{\ba}$, we can define a map
\begin{gather} 
\mathcal{P}\rightarrow\mathcal{M}_{\ba},\, p\mapsto (m_{\bd}(p))_{\bd\in\Lambda_{\ba}} \label{eq:pnmid}
\end{gather}
where $(m_{\bd}(p))_{\bd\in\Lambda_{\ba}}$ is defined for $p\in\mathcal{P}$ of the form (\ref{eq:pnmelement}) by
\begin{gather}
m_{\bd}(p):=\begin{cases} 1 & \bd=\bd(I_a,1),\, \bd(J_b,\#J_b-1),\, \bd(K_c,r_c) \\ & \textrm{ for some }1\leq a\leq A,\, 1\leq b\leq B,\, 1\leq c\leq C \\ n-r(p) & \bd=\bd(\emptyset,1) \\ 0 & \textrm{otherwise}, \end{cases} \notag
\end{gather}
and the map (\ref{eq:pnmid}) is bijective. 
\end{prop}

Consequently, for an $m$-tuple flag $F\in G^m/P^m\cong Fl_{\ba_{n,m}}(\K^n)$, if the composition of the surjection (\ref{eq:surjabstract}) and the inverse of the bijection (\ref{eq:pnmid}) maps $\mathrm{diag}(G)\cdot F$ to $p\in\mathcal{P}$ of the form (\ref{eq:pnmelement}), then by definition of these maps, $F$ splits into indecomposable $m$-tuple flags as in (\ref{eq:splitspecific}). 
Hence $F$ is $p$-constructible. 
The converse also holds. 
Hence $p\in\mathcal{P}$ is the unique element such that $F$ is $p$-constructible, and it does not depend on the choice of representatives of the orbit $\mathrm{diag}(G)\cdot F$. 
Furthermore, the induced map $\pi$ defined in (\ref{eq:mainsurj}) coincides with the composition of the surjection (\ref{eq:surjabstract}) and the inverse of (\ref{eq:pnmid}), hence it is surjective, which is the consequence of Theorem \ref{thm:main}. 

Now, each fibre $\pi^{-1}(p)$ for $p\in \mathcal{P}$ of the form (\ref{eq:pnmelement}) is identified with the fibre of the surjection (\ref{eq:surjabstract}) at $(m_{\bd}(p))_{\bd\in\Lambda_{\ba}}\in\mathcal{M}_{\ba}$. 
From (\ref{eq:preimage}), we have a homeomorphism
\begin{equation}
\pi^{-1}(p)\cong \prod_{c=1}^C\mathrm{diag}(\GL{r_c})\backslash Fl_{\bd(K_c,r_c)}^{\mathrm{ind}}(\K^{r_c}). \label{eq:fibreid}
\end{equation}
The map $\iota_p\colon\prod_{c=1}^C(\pr{r_c-1}\K)^{\#K_c-r_c-1}\rightarrow \pi^{-1}(p)$ defined in (\ref{eq:emb}) coincides with the composition of the open dense embedding \[\prod_{c=1}^C\iota_{\bd(K_c,r_c)}\colon\prod_{c=1}^C(\pr{r_c-1}\K)^{\#K_c-r_c-1}\hookrightarrow \prod_{c=1}^C \mathrm{diag}(\GL{r_c})\backslash Fl_{\bd(K_c,r_c)}^{\mathrm{ind}}(\K^{r_c})\] where $\iota_{\bd(K_c,r_c)}$ is the open dense embedding deined in (\ref{eq:embspecific}), and the homeomorphism (\ref{eq:fibreid}). 
Remark that if $C=0$, then both-hands sides of the homeomorphism (\ref{eq:fibreid}) are singletons.  
Consequently, we obtain Theorem \ref{thm:rep}. 


\section{Stabilisers and applications}

In this section, we set $G=\GL{n}$ with $n\geq 2$, and $P$ to be its maximal parabolic subgroup such that $G/P\cong\pr{n-1}\K$. 
Fix positive integers $n$ and $m$. 
Using the surjection $\pi\colon\mathrm{diag}(G)\backslash G^m/P^m\twoheadrightarrow\mathcal{P}=\mathcal{P}_{n,m}$ in Theorem \ref{thm:main}, we can observe each orbits systematically with the notion of the finite set $\mathcal{P}$. 
For instance, we can determine the conjugacy class of the stabilisers of each orbit explicitly, which depends only on $\mathcal{P}$. 
For $J\subset \{1,2,\ldots,m\}$ satisfying $3\leq\#J=:r+1$, then $g\in\GL{r}$ fixes $D(J)\in Fl_{\bd(J,\#J-1)}^{\mathrm{ind}}(\K^r)$ if and only if $g$ is a scalar matrix. 
Similarly, for $K\subset\{1,2,\ldots,m\}$ and $r\in\N$ satisfying $4\leq r+2\leq \#K$, we can take a representative of a $\mathrm{diag}(\GL{r})$-orbit in $Fl_{\bd(K,r)}^{\mathrm{ind}}(\K^r)$ such that its stabiliser coincides with the subgroup of $\GL{r}$ consisting of scalar matrices. 
Hence, we have following corollary. 
\begin{cor}
Let $n,m,G,P$, and $\mathcal{P}=\mathcal{P}_{n,m},\,\pi$ be as above. 
For $p\in\mathcal{P}$ of the form (\ref{eq:pnmelement}), 
the conjugacy class of stabilisers and dimensions of orbits do not depend on the choice of $O\in\pi^{-1}(p)$: 
we can take a representative $F\in O$ such that the stabiliser of $F$ is
\begin{equation*}
\left\{\left.
\scalebox{0.7}{$\displaystyle
\left(\begin{array}{ccccccc|c}
D &&&&&&& \multirow{7}{*}{$\scalebox{2}{$\displaystyle\ast$}$} \\
& E_1 &&&&\scalebox{1.5}{$\displaystyle0$}&& \\
&& \ddots &&&&& \\
&&& E_B &&&& \\
&&&& F_1 &&& \\
&\scalebox{1.5}{$\displaystyle0$}&&&& \ddots && \\
&&&&&& F_C & \\ \hline
\multicolumn{7}{c|}{\scalebox{1.5}{$\displaystyle0$}} & \scalebox{2}{$\displaystyle\ast$}
\end{array}\right)
$}
\right| \begin{array}{l} D=\mathrm{diag}\left(d_1,d_2,\ldots,d_A\right)\in \GL{A} \\ E_b=e_b\bm{1}\in \GL{\#J_b-1},\;e_b\in\K^{\times} \\ F_c=f_c\bm{1}\in \GL{r_b},\;f_c\in\K^{\times} \\ \textrm{for }1\leq b\leq B,\;1\leq c\leq C  \end{array}\right\}, 
\end{equation*}
and we have $\dim O=nr(p)-A-B-C$.
\end{cor}

Using this formula, we can calculate the dimensions of orbits and obtain the following corollary on the existence of open orbits. 
\begin{cor} \label{thm:open}
There exists an open $\mathrm{diag}(G)$-orbit on $G^m/P^m$ if and only if $n+1\geq m$, and open orbits for the cases where $n+1\geq m$ are given explicitly as below: 
\begin{enumerate}
\item If $n\geq m$, then $p_o=\left(\left\{ \{j\} \right\}_{j=1}^m, \emptyset ,\emptyset \right)$ is an element of $\mathcal{P}$, and $\pi^{-1}(p_o)\subset \mathrm{diag}(G)\backslash G^m/P^m$ is a singleton. 
The only orbit 
\[\mathrm{diag}(G)\cdot\left(F(\emptyset)^{\oplus(n-m)}\oplus\bigoplus_{j=1}^mF(\{j\})\right)=\mathrm{diag}(G)\cdot\left(\K^n;\K e_1, \K e_2, \ldots ,\K e_m\right)\]
contained in $\pi^{-1}(p_o)$ is an open orbit, 
\item If $n+1= m$, then $p_o=\left(\emptyset, \{\{1,2,\ldots,m\}\},\emptyset \right)$ is an element of $\mathcal{P}$, and $\pi^{-1}(p_o)\subset \mathrm{diag}(G)\backslash G^m/P^m$ is a singleton. 
The only orbit 
\[\mathrm{diag}(G)\cdot D(\{1,2,\ldots,m\})=\mathrm{diag}(G)\cdot\left(\K^n;\K e_1, \K e_2, \ldots ,\K e_n, \K\left(\sum_{i=1}^ne_i\right)\right) \]
contained in $\pi^{-1}(p_o)$ is an open orbit, 
\end{enumerate}
where $\mathcal{P}, \pi, F(I)$, and $D(J)$ are those in Theorem \ref{thm:rep}. 
\end{cor}
Not only the existence of open orbits, but we can also check the infiniteness of orbits easily. 
As we have seen in Theorem \ref{thm:rep}, for the surjection $\pi$ from $\mathrm{diag}(G)\backslash G^m/P^m$ to the finite set $\mathcal{P}$ and an element $p\in\mathcal{P}$ of the form (\ref{eq:pnmelement}), the fibre $\pi^{-1}(p)$ is a singleton if the $3$rd family of $p$ is empty, and if it is not empty, then the fibre $\pi^{-1}(p)$ has uncoutably many elements. 
Hence, finiteness of $\mathrm{diag}(G)\backslash G^m/P^m$ is equivalent to the bijectivity of $\pi$, and it occurs if and only if the $3$rd families of all elements of $\mathcal{P}$ are empty. 
From Definition \ref{thm:pnm}, we can see that it is equivalent to $m\leq 3$, and have the following corollary: 
\begin{cor} \label{thm:finite}
It is equaivalent that $m\leq 3$ and there are only finitely many $\mathrm{diag}(G)$-orbits on $G^m/P^m$. 
\end{cor}
This result corresponds to the finiteness of isomorphism classes of indecomposable representations on quivers of Dynkin type. 
In our setting, orbits on the multiple flag variety $G^m/P^m$ can be identified with isomorphism classes of representations of the quiver of type $A_2$ if $m=1$, $A_3$ if $m=2$, and $D_4$ if $m=3$. 
If $m=4$, then the corresponding quiver to $G^m/P^m$ is the extended Dynkin diagram of type $\hat{D}_4$, which is a tame quiver. 
In this case, from Theorem \ref{thm:rep}, each fibre of $\pi$ satisfies either it is a singleton or it has an open dense embedding from $\pr{1}\K$.
It corresponds to the fact that isomorphism classes of indecomposable representations in every dimension vector on a tame quiver are described by finite number of one-parameter families. 
If $m\geq 5$, then the corresponding quiver to $G^m/P^m$ is wild. 
In these cases, there exists an element $p\in\mathcal{P}$ whose fibre $\pi^{-1}(p)$ has an open dense embedding from a direct product of $\pr{1}\K$ of at least twice, which is an at least two-parameters family. 



\end{document}